\documentclass[11pt]{article}
\usepackage{latexsym}
\usepackage{amsfonts}
\usepackage{graphicx}
\usepackage{epstopdf}
\usepackage{xcolor}
\usepackage[sort]{natbib}
\usepackage{lineno}

\newtheorem{theorem}{Theorem}
\newtheorem{lemma}{Lemma}
\newtheorem{corol}{Corollary}
\newtheorem{property}{Property}

\newcommand{\be}{\begin{equation}}
\newcommand{\ee}{\end{equation}}

\usepackage[normalem]{ulem}

\begin{document}
\title{\bf Partitioning Items Into Mutually Exclusive Groups}
\author{Pawel Kalczynski, Zvi Goldstein, and Zvi Drezner\\
Steven G. Mihaylo College of Business and Economics\\
California State University-Fullerton\\
Fullerton, CA 92834.\\e-mail: pkalczynski@fullerton.edu; zgoldstein@fullerton.edu; zdrezner@fullerton.edu}
\date{}
\maketitle

\begin{abstract}
We investigate a new model for partitioning a set of items into groups (clusters). The number of groups is given and the distances between items are well defined. These distances may include weights. The sum of the distances between all members of the same group is calculated for each group, and the objective is to find the partition of the set of items that minimizes the sum of these individual sums. Two problems are formulated and solved. In the first problem the number of items in each group are given. For example, all groups must have the same number of items. In the second problem there is no restriction on the number of items in each group.

We propose an optimal algorithm for each of the two problems as well as a heuristic algorithm. Problems with up to 100 items and 50 groups are tested. In the majority of instances the optimal solution was found using IBM's CPLEX. The heuristic approach, which is very fast, found all confirmed optimal solutions and equal or better solutions when CPLEX was stopped after five hours. The problem with given group sizes can also be formulated and solved as a quadratic assignment problem.
\end{abstract}

\renewcommand{\baselinestretch}{1.6}
\renewcommand{\arraystretch}{0.625}
\large
\normalsize

\section{Introduction}

Consider sport teams, such as those competing in the World cup, Olympics, NBA, NFL, NHL, various college teams. The teams need to be partitioned into divisions. In the first round teams compete within the division. In the second round the winners (and in some cases also runners-up) advance to compete in the finals. We are interested in a good partition of the teams into divisions for the first round of competition. Each team competes with each other team in the same division at least twice: once at home and once at the home of the other team. We are therefore interested to form divisions so that the total travel time by the teams is minimized.
The number of divisions is pre-determined.	
It may be required to have a given number of teams in each division or the number of teams in each division is flexible.	

This problem is similar to many multi-facility location problems. In facility location problems, a facility is located at the ``center" of each division and the service distances are the distances to the central facility of each division (group) rather than pair-wise distances without a central facility. This may be the appropriate model if the competition within the division is held at a central location.

In typical multi-facility location models, a given number of facilities, $p$, need to be located. For each demand point, its distance is defined as the distance to the closest facility. The following is a list of common objectives. 

\begin{enumerate}
\item The objective of the $p$-median problem \citep{Das95,KH79med,DM15}, sometimes refereed to as the multi-source Weber problem \citep{BHMT00,KS72}, is to minimize the sum of the weighted distances to all demand points. The 1-median problem is termed the Weber problem \citep{W09} has a long history dating back to Fermat in the 1600s \citep{W93,DKSW02,C19}. \citet{C19} has references dating back to 1687.

\item The $p$-median problem with equal weights and squared Euclidean distances is used in many papers to model clustering problems \citep[e.g., ][]{L57,HW79,M67a,A09}. For recent reviews see \cite{BODX15,PABM19}.
	
\item The objective of the $p$-center problem \citep{KH79cen,Pap25,CLY15} is to minimize the largest distance to all demand points. The 1-center problem is termed the minimax location problem \citep{Sylv,Sylv2,EH72,Dr11}.

\item One of the competitive location models is based on the assumption that customers patronize the closest facility \citep{H29,Dr82,Hak83}. The objective is to locate one or more facilities to attract the maximum buying power from customers located at the demand points.

\item Obnoxious facility models \citep{CG78,ER90,Erk:89,DW96} assume that facilities generate nuisance and facilities should be located as far as possible from demand points. In most models it is assumed that the nuisance affecting a demand point is generated by the closest facility

\item The objective of equitable location models \citep{BBKW07,BDTW09} is to equalize as much as possible the number of demand points attracted to each facility.

\item A wide variety of equity objectives were proposed. \citet{Ma86} suggested to minimize the variance of distances, \citet{DrThErk} suggested to minimize the range of distances, \citet{DDG09a} suggested to minimize the Gini Coefficient of the Lorenz curve, \citet{Pap257} suggested to minimize the Quintile share ratio. \citet{EL95} list nineteen equity measures. As far as we know no one considered the multiple facility extension of these problems.

\end{enumerate}

In all these multiple facilities location models, for a given set of locations for the facilities, the demand points are partitioned into groups. The plane is partitioned into polygons which is termed a Voronoi diagram \citep{SI92,V08,SOK95,OKSC00}. Each facility is surrounded by a polygon and provides its services to the demand points located in that polygon. The set of demand points is partitioned into groups.

In this paper we propose a new criterion for grouping the set of demand points. The objective is to minimize the total intra-group distances. This means, for each group the total distance between all pairs of demand points belonging to the same group is calculated and the sum of these total distances is minimized. The number of groups is pre-specified. Clearly, if the number of groups increases, the objective function decreases. In fact, when the number of groups is equal to the number of demand points, each demand point belongs to its own group, and the objective function is zero.

Weights can be easily incorporated into the models. If the population at each demand point is a factor, the distance between two demand points can be multiplied by, for example, the product or sum of the population counts at the two demand points thus defining a new ``distance" value.

Two models are proposed. In Model A the number of items in each group is pre-specified. If there are, for example, 12 items to be partitioned into 3 groups, we can require that one group will contain 5 items, one 4 items, and one 3 items. It is possible that all groups are required to have the same number of items. An equal division of the items in each group may be required when sport teams are partitioned into divisions and we wish to have the same number of teams in each division. In Model B, the number of items in each group is not specified and each group can include any number of items as long as every item belongs exclusively to one group.

The paper is organized as follows. In Section \ref{sec2} the fixed group size model is formulated. The quadratic assignment problem which can be used to solve this problem as well as a non-trivial binary linear programming formulation are detailed. In Section \ref{sec3} the variable size model is presented and a binary linear program formulated. The main contribution of this paper is the effective heuristic procedures which are detailed in Section \ref{sec4}. Two case studies are described in Section \ref{sec5}, and extensive computational experiments in Section \ref{sec6} demonstrate the effectiveness of the heuristic approaches. We conclude the paper in Section \ref{sec7} suggesting the application of the proposed effective heuristic approach to other multi-facility location problems as well.

\section{\label{sec2}Model A: Fixed Group Sizes}

There are $n$ items. The distance between items $i\ne j$ is $d_{ij}$. By definition $d_{ii}=0$. We need to create $p$ groups. The number of items in these groups are required to be $n_1,n_2,\ldots,n_p$ such that $\sum\limits_{i=1}^pn_i=n$. For generality we do not assume that each group has the same number of items. It is possible that we allow a certain number of items (can be any of the $n$ items) to be excluded from the groups. In this case $\sum\limits_{i=1}^pn_i<n$.

When all groups are required to be of the same size $\frac n p$, the total number of combinations is by the multinomial distribution: $\frac{n!}{\left(\frac np !\right)^pp!}$ which can be quite small for small values of $n$ and $p$, allowing for total enumeration. For example, for $n=12$, $p=2$ there are 462 possible partitions of items.

\subsection{Models Related to Model A}

Model A can be formulated as a quadratic assignment problem \citep[QAP,][]{Dr15,Gil62,Law,Rendl,KB57}, and has similar features to the gray pattern problem \citep{Ta95} which is a special case of the quadratic assignment problem.

\subsubsection{The Quadratic Assignment Problem (QAP)}

The quadratic assignment problem  is considered to be one of the most difficult optimization problems to solve optimally.
It is defined as follows.
A set of $n$ possible sites are given and $n$ facilities are to be located on these sites, one facility at a site. Let $c_{ij}$ be the cost per unit distance between facilities $i$ and $j$ and  $d_{ij}$ be the distance between sites $i$ and $j$. The cost $f$ to be minimized over all possible permutations, calculated for an assignment of facility $i$ to site $p(i)$ for $i=1,\ldots,n$, is:
\begin{equation}\label{eq1}
f=\sum\limits_{i=1}^n\sum\limits_{j=1}^n c_{ij}d_{p(i)p(j)}
\end{equation}

\subsubsection{The Gray Pattern Problem}

In the gray pattern problem \citep{Ta95}, the objective is to select $k$ out of $n$ black dots so that the pattern will be as uniform as possible. This is achieved by a specially designed distances matrix. 

\citet{Ta95} suggested the following cost matrix $\{c_{ij}\}$: $c_{ij}=1$ For $1\le i,j\le k$ and all other $c_{ij}=0$. Such a matrix in (\ref{eq1}) results in the sum of distances between all pairs of facilities for $1\le i\ne j\le k$ because $d_{ii}=0$. The solution to the resulting QAP is the best possible selection of $k$ out of the $n$ items. The cost matrix $\{c_{ij}\}$ can be visualized as the top left corner of the cost matrix has a $k$ by $k$ square of 1's, and the rest of the cost matrix is all zeros.

The distances of the gray pattern problem suggested by \cite{Ta95} are based on a rectangle of dimensions $n_1$ by $n_2$. A gray pattern of $p$ black points is selected from the $n=n_1\times n_2$ points in the rectangle while the rest of the points remain white. The objective is to have a gray pattern where the black points are distributed as uniformly as possible. To achieve that, \citet{Ta95} created a distance matrix $d_{ij}$ between locations $i$ and $j$ by minimum entropy in Physics.
Let $1\le i\le n$ be the serial number of a point in the rectangle. Define $r(i)$ the row of point $i$, and $c(i)$ the column of point $i$:
\begin{equation}\label{rs}
r(i)=\lfloor\frac{i-1}{n_2}\rfloor;~c(i)=i-1-n_2r(i).
\end{equation}
The distance between points $i$ and $j$ ($i\ne j$) is
\begin{equation}\label{dis}
d_{ij}=\max\limits_{u,v\in\{-1,0,1\}}\left\{\frac{100,000}{(r(i)-r(j)+un_1)^2+(c(i)-c(j)+vn_2)^2}\right\}
\end{equation}
The distances are rounded to the nearest integer.  It is interesting to note that when the result of (\ref{dis}) is exactly an integer+$\frac12$ (such as 1562.5) the distance is rounded down which is not the customary rounding rule. We also define $d_{ii}=0$ (it is infinite by (\ref{dis})).

\citet{DMP14} optimally solved gray pattern problems of $n=64$ (8$\times8$) for $2\le p\le 25$, $n=256$ (16$\times16$) for $2\le p\le 12$,  $n=576$ (24$\times24$) for $2\le p\le 9$, and  $n=1024$ (32$\times32$) for $2\le p\le 8$.

The gray pattern problem can be extended to patterns of several colors, not only black points. The sequence of several squares along the diagonal can be used to create a pattern of several colors with certain percentage for each color so that the pattern will be as uniform as possible. Points that are not selected are white like in the gray pattern problem. This idea is used to formulate Model A as a quadratic assignment problem.

\subsection{Formulating Model A as a Quadratic Assignment Problem}

The quadratic assignment formulation for the grouping problem is based on the gray pattern model \citep{Ta95,Dr06,DK17,DMP14}. 
\citet{Dr06} proposed the one group problem:  selecting the best one group of $k$ items out of $n$ available items and observed that this problem is equivalent to the gray pattern problem.

For our problem of multiple groups we generate a cost matrix $\{c_{ij}\}$ in Equation (\ref{eq1}) that consists of $p$ ``squares" of 1's on the diagonal of the matrix of sizes $n_1, n_2,\ldots,n_p$ and the rest of the $c_{ij}=0$. This way, only distances between pairs of items belonging to the same group are counted in the objective function. If $\sum\limits_{i=1}^pn_i=n$, the whole diagonal is covered with squares and if $\sum\limits_{i=1}^pn_i<n$, like the case of finding one group \citep{Dr06}, some items are not selected for any group.

For example, if 12 items are split into three groups of 4 items each, the $C=\{c_{ij}\}$ matrix is depicted in Figure \ref{12}. The distance matrix consists of the distances between items and not defined by (\ref{dis}).

\begin{figure}[ht!]
$$
\left(\begin{array}{c}
1~~ 1 ~~1 ~~1 ~~0~~ 0~~0~~ 0 ~~0 ~~0 ~~0 ~~0\\
1~~ 1 ~~1 ~~1 ~~0~~ 0~~0~~ 0 ~~0 ~~0 ~~0 ~~0\\
1~~ 1 ~~1 ~~1 ~~0~~ 0~~0~~ 0 ~~0 ~~0 ~~0 ~~0\\
1~~ 1 ~~1 ~~1 ~~0~~ 0~~0~~ 0 ~~0 ~~0 ~~0 ~~0\\
0~~ 0 ~~0 ~~0 ~~1~~ 1~~1~~ 1 ~~0 ~~0 ~~0 ~~0\\
0~~ 0 ~~0 ~~0 ~~1~~ 1~~1~~ 1 ~~0 ~~0 ~~0 ~~0\\
0~~ 0 ~~0 ~~0 ~~1~~ 1~~1~~ 1 ~~0 ~~0 ~~0 ~~0\\
0~~ 0 ~~0 ~~0 ~~1~~ 1~~1~~ 1 ~~0 ~~0 ~~0 ~~0\\
0~~ 0 ~~0 ~~0 ~~0~~ 0~~0~~ 0 ~~1 ~~1 ~~1 ~~1\\
0~~ 0 ~~0 ~~0 ~~0~~ 0~~0~~ 0 ~~1 ~~1 ~~1 ~~1\\
0~~ 0 ~~0 ~~0 ~~0~~ 0~~0~~ 0 ~~1 ~~1 ~~1 ~~1\\
0~~ 0 ~~0 ~~0 ~~0~~ 0~~0~~ 0 ~~1 ~~1 ~~1 ~~1\\
\end{array}
\right)$$

\caption{\label{12}The $C$ matrix for $n=12$; $p=3$; and 4 items in each group}	
\end{figure}

\subsection{Binary Linear Programming Formulation for Model A}

The grouping model can also be formulated as a binary linear program. Define $np$ binary variables $X_{ik}$ such that $X_{ik}=1$ if  item $i$ belongs to group $k$ and 0 otherwise. In addition we define $n^2$ variables $Y_{ij}$ such that $Y_{ij}=1$ if items $i$ and $j$ belong to the same group and 0 otherwise. The number of the $Y$ variables can be reduced from $n^2$ to $\frac12n(n-1)$ by observing that by definition $Y_{ii}=1$ and $Y_{ij}=Y_{ji}$.

\begin{eqnarray}
&&\min\left\{\sum\limits_{i=1}^n\sum\limits_{j=1}^nd_{ij}Y_{ij}\right\}\label{eq2}\\
\mbox{Subject to:}&&\nonumber\\
&&-M(1-Y_{ij})\le \sum\limits_{k=1}^pkX_{ik}-\sum\limits_{k=1}^pkX_{jk}\le M(1-Y_{ij})\mbox{ for }i,j=1,\ldots,n\label{eq3}\\
&&\sum\limits_{i=1}^n\sum\limits_{j=1}^nY_{ij}=\sum\limits_{k=1}^p n_k^2\label{eq4}\\
&&\sum\limits_{k=1}^pX_{ik}=1 \mbox{ for }1\le i\le n\label{eq5}\\
&&\sum\limits_{i=1}^nX_{ik}=n_k\mbox{ for }1\le k\le p\label{eq6}\\
&&Y_{ij},X_{ik}\in\{0,1\}.
\end{eqnarray}
 
Equation (\ref{eq2}) guarantees that only distances between items belonging to the same group are counted in the objective function. The center term of inequality (\ref{eq3}) is equal to 0 if items $i$ and $j$ belong to the same group, and is positive or negative otherwise. Its maximum possible absolute value is $p-1$. We use $M=p$, so that if $Y_{ij}=0$ the constraint is always satisfied and if $Y_{ij}=1$ it is only satisfied when items $i$ and $j$ belong to the same group. Equation (\ref{eq4}) forces  $Y_{ij}=1$ for items belonging to the same group because there are exactly $\sum\limits_{k=1}^p n_k^2$ such pairs. Equation (\ref{eq5}) guarantees that each item belongs to one and only one group, and equation (\ref{eq6}) guarantees that group $k$ has $n_k$ items.
	
\subsubsection{Reducing the Number of Variables}
	
Since $Y_{ij}=Y_{ji}$ and $Y_{ii}=1$ the problem can be reduced in size by having as variables $Y_{ij}$ for $j>i$. We get:

\begin{eqnarray}
&&\min\left\{2\sum\limits_{i=1}^{n-1}\sum\limits_{j=i+1}^nd_{ij}Y_{ij}\right\}\nonumber\\
\mbox{Subject to:}&&\label{formul}\\
&&-M(1-Y_{ij})\le \sum\limits_{k=1}^pkX_{ik}-\sum\limits_{k=1}^pkX_{jk}\le M(1-Y_{ij})\mbox{ for }j>i\nonumber\\
&&\sum\limits_{i=1}^{n-1}\sum\limits_{j=i+1}^nY_{ij}=\sum\limits_{k=1}^p \frac12n_k(n_{k}-1)\nonumber\\
&&\sum\limits_{k=1}^pX_{ik}=1 \mbox{ for }1\le i\le n\nonumber\\
&&\sum\limits_{i=1}^nX_{ik}=n_k\mbox{ for }1\le k\le p\nonumber\\
&&Y_{ij}\mbox{ for }j>i,X_{ik}\in\{0,1\}.\nonumber
\end{eqnarray}

In the computational experiments we used this more compact formulation.

\section{\label{sec3}Model B: Variable Group Sizes}
Model B cannot be formulated as a quadratic assignment problem because the structure of the $C$ matrix is not determined. The BLP formulation is more complicated because the numbers $n_k$ are variables rather than having given values.

The total number of possible groups $P(n,p)$, i.e., the number of ways of partitioning a set of $n$ elements into $p$ non-empty subsets, satisfies the recursion: $P(n,p)=P(n-1,p-1)+pP(n-1,p)$. Suppose that item $n$ is added to the groups formed by the first $n-1$ items. The first term $P(n-1,p-1)$ is the number of possible groups when item $n$ forms a new group of one item. The second term $pP(n-1,p)$ is the number of possible groups where item $n$ joins one of the existing groups. $P(n,p)$ can be easily calculated recursively, for example, by Microsoft Excel, starting with $P(n,1)=P(p,p)=1$.

This recursion leads to the Stirling number of the second kind \citep{S1764,AS72}: $$P(n,p)=\frac{1}{p!}\sum\limits_{j=1}^p (-1)^{p-j}\frac{p!}{j!(p-j)!}j^n$$ which can be quite small for small values of $n$ and $p$, allowing for total enumeration. For example, for $n=12$, $p=2$ there are 2,047 possible partitions of items. Note that for $p=2$, $P(n,2)=2^{n-1}-1$.

\subsection{Binary Linear Programming Formulation for Model B}

In Model B, $n_k$ is a variable that can assume a value between 1 and $n-p+1$. Each group must have at least one item and $\sum\limits_{k=1}^pn_k=n$. Similarly to the Model A formulation (\ref{formul}), Model B can be formulated as a BLP by adding $(n-p+1)p$ binary variables $U_{ik}$ for $i=1,\ldots n-p+1; k=1,\ldots p$. $U_{ik}=1$ if $i=n_k$ ($n_k$ is a result of the solution) and 0 otherwise. Consequently, $n_k=\sum\limits_{i=1}^{n-p+1}iU_{ik}$.

\begin{eqnarray}
&&\min\left\{\sum\limits_{i=1}^{n-1}\sum\limits_{j=i+1}^nd_{ij}Y_{ij}\right\}\\
\mbox{Subject to:}&&\nonumber\\
&&-M(1-Y_{ij})\le \sum\limits_{k=1}^pkX_{ik}-\sum\limits_{k=1}^pkX_{jk}\le M(1-Y_{ij})\mbox{ for }j>i\\
&& n+2~\sum\limits_{i=1}^{n-1}\sum\limits_{j=i+1}^nY_{ij}=\sum\limits_{i=1}^{n-p+1} i^2~\sum\limits_{k=1}^p U_{ik}\\
&&\sum\limits_{k=1}^pX_{ik}=1 \mbox{ for }1\le i\le n\\
&&\sum\limits_{i=1}^nX_{ik}=\sum\limits_{i=1}^{n-p+1} i~ U_{ik}\mbox{ for }1\le k\le p\\
&&\sum\limits_{i=1}^{n-p+1} U_{ik}=1\mbox{ for }1\le k\le p\\
&&Y_{ij}\mbox{ for }j>i,X_{ik},U_{ik}\in\{0,1\}.
\end{eqnarray}

\section{\label{sec4}Heuristic Approaches}

The following heuristic approach is simpler for Model B so we present it first for Model B and then show the modifications required in order to apply it for the solution of Model A. 

\subsection{Heuristic Approach to Model B}

\begin{description}
	
	\item[Phase 1 (selecting one item for each group):]
	~
	\begin{itemize}
		\item Select the two items that are farthest from one another to be in group \#1 and group \#2. It is unlikely that the farthest two items will be in the same group.
		\item Select for the next group the item that its minimum distance to the already selected items is the largest.
		\item Continue until $p$ items are assigned to $p$ groups. 
	\end{itemize}

	\item[Phase 2 (adding items to groups):]
	~
	\begin{itemize}
		\item Check all unassigned items to be added to each of the groups. Select the item that adding it to one of the groups will increase the objective function the least and add this item to the particular group.
		\item Keep adding items to groups until all the items are assigned a group.
	\end{itemize}
	
	\item[Phase 3 (a descent algorithm):]
	~
	\begin{enumerate}
		\item \label{st31} Evaluate all combinations of moving an item from one group to another.
		\item If an improvement is found, perform the best move and go to Step \ref{st31}.
		\item If no improvement is found, stop.
	\end{enumerate}

\item[Applying GRASP:]
~

The algorithm described above yields one final solution. We introduce randomness into the procedure so that the algorithm can be repeated many times in a multi-start heuristic approach. The idea follows the ``Greedy Randomized Adaptive Search Procedure" (GRASP) suggested by \citet{FR95}. In each step in the first two phases, we find the best move and the second best move. The best move is selected with probability $\frac23$ and the second best with probability $\frac13$.
	
\end{description}

\subsection{Modification to the Heuristic Approach for Solving Model A}

We detail the modifications for equal size $s$ for all groups. It can be further modified to a list of given sizes which are not necessarily equal.

Phase 1 is not altered. In Phase 2: once a group reaches size $s$ it is no longer considered for adding more items to it. In Phase 3: since the groups must retain the same number of items, rather than evaluating the move of an item from one group to another, we evaluate all pair exchanges between groups. The GRASP approach is not altered.

\subsection{Efficient Calculations}

In order to calculate the increase in the value of the objective function when an item is added to a group, we just need to calculate the sum of the distances between the added item and all the existing items in the group.

In order to calculate the decrease  in the value of the objective function when an item is removed from a group (needed for the decent phase), we just find the sum of distances between the removed item and all the remaining items in the group.

\section{\label{sec5}Case Studies}

\subsection{Creating Divisions for Teams}
\citet{WA01} used the example of creating divisions for NBA teams to illustrate the evolutionary procedure available  in Microsoft Excel's Solver.
There are $n$ competing sport teams that play against one another during the regular season. The teams are partitioned into $p$ divisions with $n_i$ teams in each division. Each pair of teams in a division play two games against one another, one in each home team city. The problem is to form divisions in such a way that the total distance traveled by the teams is minimized.

\subsection{Dividing the Largest US Cities into Groups}

\begin{figure}[htp!]
	\begin{center}
		\setlength{\unitlength}{1in}
		\begin{picture}(4.5,3)
		\includegraphics[width=4.5in]{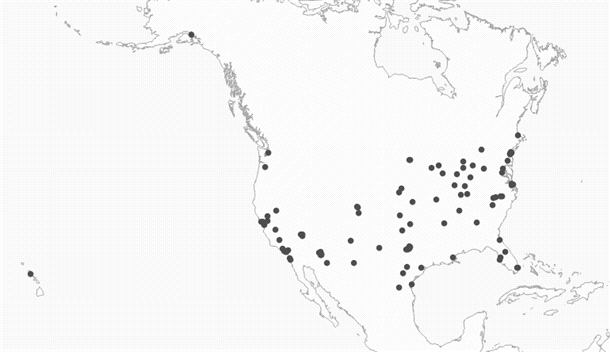}
		\end{picture}
		\vspace{0.1in}
		\caption{\label{100}The 100 Largest U.S. Cities}
	\end{center}
\end{figure}

The website {\tt https://simplemaps.com/data/us-cities} has data for more than 36,000 U.S. municipalities. The data provide both population and population proper (i.e., not including the population in the metropolitan area). We selected the proper population and therefore the largest city is New York, the fourth largest is Brooklyn, the fifth largest is Queens, and so on. We extracted data for latitude, longitude and population of the largest 100 U.S. cities. The data for these cities sorted by their population size is given in Table \ref{cities} and depicted in Figure \ref{100}. It can be used to solve problems for any number of cities up to 100. Instances based on fewer than 100 cities, use the largest ones listed in Table \ref{cities}.

The distances between any two cities can be calculated by the ``great circle" formula given in \cite{Pap3}. The distance $d$ between two cities whose latitudes are $\phi_1, \phi_2$ and longitudes $\theta_1, \theta_2$ is:
\begin{equation}\label{dist}
d=R\arccos\left(\cos\phi_1\cos\phi_2\cos(\theta_1-\theta_2)+\sin\phi_1\sin\phi_2\right)
\end{equation}
where $R=3959$ miles is the earth's radius. This formula can be rewritten to avoid large round-off errors for small distances when $\cos\frac{d}{R}\approx 1$. The  identity
$
\cos\alpha=1-2\sin^2\frac{\alpha}{2}
$ is used.
We get
\begin{eqnarray}
&\cos\phi_1\cos\phi_2\cos(\theta_1-\theta_2)+\sin\phi_1\sin\phi_2=\cos\phi_1\cos\phi_2+\sin\phi_1\sin\phi_2\nonumber\\& -(1-\cos(\theta_1-\theta_2))
\cos\phi_1\cos\phi_2= \cos(\phi_1-\phi_2)-2\sin^2\frac{\theta_1-\theta_2}{2}\cos\phi_1\cos\phi_2\nonumber
\end{eqnarray}
yielding
\begin{eqnarray}
\cos\frac{d}{R}&=&1-2\sin^2\frac{d}{2R}= \cos(\phi_1-\phi_2)-2\sin^2\frac{\theta_1-\theta_2}{2}\cos\phi_1\cos\phi_2\nonumber\\
&\rightarrow&2\sin^2\frac{d}{2R}=2\sin^2\frac{\phi_1-\phi_2}{2}+2\sin^2\frac{\theta_1-\theta_2}{2}\cos\phi_1\cos\phi_2\nonumber\\
&\rightarrow&d=2R\arcsin\sqrt{\sin^2\frac{\phi_1-\phi_2}{2}+\sin^2\frac{\theta_1-\theta_2}{2}\cos\phi_1\cos\phi_2}\label{dist1}
\end{eqnarray}

\begin{table}[ht!]
	\begin{center}
		\caption{\label{cities} List of 100 most populated U.S. Cities}
		\medskip
\small
\setlength{\tabcolsep}{1.5pt}
		\begin{tabular}{|l|r|r||l|r|r||l|r|r|}
			\hline
			City&Lat.&Long.&	City&Lat.&Long.&	City&Lat.&Long.\\
				\hline
New York&	40.6943&	-73.9249&	Milwaukee&	43.0640&	-87.9669&	Saint Paul&	44.9477&	-93.1040\\
Los Angeles&	34.1140&	-118.4068&	Albuquerque&	35.1055&	-106.6476&	Cincinnati&	39.1412&	-84.5060\\
\cline{7-9}
Chicago&	41.8373&	-87.6861&	Tucson&	32.1558&	-110.8777&	Anchorage&	61.1508&	-149.1091\\
Brooklyn&	40.6501&	-73.9496&	Fresno&	36.7834&	-119.7933&	Henderson&	36.0145&	-115.0362\\
Queens&	40.7498&	-73.7976&	Sacramento&	38.5666&	-121.4683&	Greensboro&	36.0960&	-79.8275\\
Houston&	29.7871&	-95.3936&	Mesa&	33.4016&	-111.7180&	Plano&	33.0502&	-96.7487\\
\cline{4-6}
Manhattan&	40.7834&	-73.9662&	Kansas City&	39.1239&	-94.5541&	Newark&	40.7242&	-74.1724\\
Phoenix&	33.5722&	-112.0891&	Staten Island&	40.5834&	-74.1496&	Lincoln&	40.8102&	-96.6808\\
Philadelphia&	40.0076&	-75.1340&	Atlanta&	33.7627&	-84.4231&	Toledo&	41.6639&	-83.5820\\
San Antonio&	29.4722&	-98.5247&	Long Beach&	33.8059&	-118.1610&	Orlando&	28.4801&	-81.3448\\
\cline{1-3}
Bronx&	40.8501&	-73.8662&	Colorado Springs&	38.8673&	-104.7605&	Chula Vista&	32.6281&	-117.0144\\
San Diego&	32.8312&	-117.1225&	Raleigh&	35.8323&	-78.6441&	Irvine&	33.6772&	-117.7738\\
\cline{7-9}
Dallas&	32.7938&	-96.7659&	Miami&	25.7840&	-80.2102&	Fort Wayne&	41.0885&	-85.1436\\
San Jose&	37.3020&	-121.8488&	Virginia Beach&	36.7335&	-76.0435&	Jersey City&	40.7161&	-74.0683\\
Austin&	30.3038&	-97.7545&	Omaha&	41.2634&	-96.0453&	Durham&	35.9801&	-78.9045\\
Jacksonville&	30.3322&	-81.6749&	Oakland&	37.7903&	-122.2165&	Saint Petersburg&	27.7931&	-82.6652\\
\cline{4-6}
San Francisco&	37.7561&	-122.4429&	Minneapolis&	44.9635&	-93.2679&	Laredo&	27.5536&	-99.4890\\
Columbus&	39.9859&	-82.9852&	Tulsa&	36.1284&	-95.9037&	Buffalo&	42.9016&	-78.8487\\
Indianapolis&	39.7771&	-86.1458&	Arlington&	32.6998&	-97.1251&	Madison&	43.0809&	-89.3921\\
Fort Worth&	32.7813&	-97.3466&	New Orleans&	30.0687&	-89.9288&	Lubbock&	33.5665&	-101.8867\\
\cline{1-3}
Charlotte&	35.2080&	-80.8308&	Wichita&	37.6894&	-97.3440&	Chandler&	33.2828&	-111.8517\\
Seattle&	47.6217&	-122.3238&	Cleveland&	41.4766&	-81.6805&	Scottsdale&	33.6872&	-111.8650\\
\cline{7-9}
Denver&	39.7621&	-104.8759&	Tampa&	27.9937&	-82.4454&	Glendale&	33.5797&	-112.2246\\
El Paso&	31.8478&	-106.4310&	Bakersfield&	35.3528&	-119.0354&	Reno&	39.5487&	-119.8486\\
Washington&	38.9047&	-77.0163&	Aurora&	39.7085&	-104.7274&	Norfolk&	36.8945&	-76.2590\\
Boston&	42.3189&	-71.0838&	Honolulu&	21.3293&	-157.8460&	Winston-Salem&	36.1029&	-80.2610\\
\cline{4-6}
Detroit&	42.3834&	-83.1024&	Anaheim&	33.8390&	-117.8572&	North Las Vegas&	36.2880&	-115.0901\\
Nashville&	36.1714&	-86.7844&	Santa Ana&	33.7366&	-117.8819&	Irving&	32.8584&	-96.9702\\
Memphis&	35.1047&	-89.9773&	Corpus Christi&	27.7173&	-97.3822&	Chesapeake&	36.6778&	-76.3024\\
Portland&	45.5372&	-122.6500&	Riverside&	33.9382&	-117.3949&	Gilbert&	33.3103&	-111.7463\\
\cline{1-3}
Oklahoma City&	35.4677&	-97.5138&	Lexington&	38.0423&	-84.4587&	Hialeah&	25.8696&	-80.3045\\
Las Vegas&	36.2288&	-115.2603&	Saint Louis&	38.6358&	-90.2451&	Garland&	32.9100&	-96.6305\\
\cline{7-9}
Louisville&	38.1662&	-85.6488&	Stockton&	37.9766&	-121.3111&	&	&	\\
Baltimore&	39.3051&	-76.6144&	Pittsburgh&	40.4396&	-79.9763&	&	&	\\
			\hline
		\end{tabular}
	\end{center}
\end{table}

\newpage

\begin{table}[ht!]
	\begin{center}
		\caption{\label{reseq} Results for fixed group sizes}
		\medskip
		\begin{tabular}{|r|r|c|r|r|r||c|c||c|c|}
			\hline
			&&Combin-&\multicolumn{3}{|c||}{Heuristic Runs}&\multicolumn{2}{|c||}{CPLEX}&\multicolumn{2}{|c|}{QAP}\\
			\cline{4-10}			
			$n$&$k$&ations$^*$&Best&Times$\dagger$&Time$\ddagger$&\%&+&Times$\dagger$&Time$\ddagger$\\
			\hline
12&	2&	462&	46,323.8&	10,000&	0.02&	0&	OPT&	10&	0.00\\
12&	3&	5,775&	23,476.8&	10,000&	0.03&	0&	OPT&	10&	0.00\\
12&	4&	15,400&	12,172.4&	10,000&	0.02&	0&	OPT&	10&	0.00\\
12&	6&	10,395&	5,009.4&	9,512&	0.02&	0&	OPT&	10&	0.00\\
18&	2&	24,310&	118,124.2&	10,000&	0.09&	0&	OPT&	10&	0.02\\
18&	3&	2,858,856&	57,393.8&	10,000&	0.14&	0&	OPT&	10&	0.02\\
18&	6&	1.91E+08&	13,965.6&	3,442&	0.09&	0&	OPT&	10&	0.02\\
18&	9&	34,459,425&	5,830.8&	10,000&	0.08&	0&	OPT&	10&	0.03\\
30&	2&	77,558,760&	336,107.2&	10,000&	0.50&	0&	OPT&	10&	0.15\\
30&	3&	9.25E+11&	161,559.2&	10,000&	0.41&	0&	OPT&	10&	0.15\\
30&	5&	1.14E+16&	67,421.0&	8,864&	0.50&	0&	OPT&	10&	0.15\\
30&	6&	1.23E+17&	50,457.8&	3,269&	0.36&	0&	OPT&	10&	0.16\\
30&	10&	1.21E+18&	18,078.4&	4,104&	0.44&	0&	UB$^*$&	9&	0.24\\
30&	15&	6.19E+15&	7,007.2&	3,333&	0.30&	0&	OPT&	8&	0.29\\
\hline									
40&	2&	6.89E+10&	501,424&	10,000&	0.61&	0&	OPT&	10&	0.42\\
40&	4&	1.96E+20&	149,708&	10,000&	1.67&	0&	OPT&	10&	0.42\\
40&	5&	6.38E+22&	102,882&	9,999&	1.20&	0&	OPT&	10&	0.43\\
40&	8&	4.71E+26&	44,976&	6,982&	1.30&	1.79\%&	UB&	10&	0.44\\
40&	10&	3.55E+27&	32,782&	7,891&	1.02&	6.71\%&	UB&	10&	0.56\\
40&	20&	3.20E+23&	7,082&	2,296&	0.83&	0&	OPT&	4&	0.84\\
50&	2&	6.32E+13&	801,378&	10,000&	0.95&	0&	OPT&	10&	1.00\\
50&	5&	4.03E+29&	182,112&	5,243&	3.47&	0.48\%&	UB&	10&	1.01\\
50&	10&	1.35E+37&	53,164&	1,110&	1.75&	0&	UB&	10&	1.08\\
50&	25&	5.84E+31&	6,782&	6,994&	1.48&	0&	OPT&	0&	1.89\\
100&	2&	5.04E+28&	3,656,540&	10,000&	11.89&	0&	OPT&	10&	26.71\\
100&	4&	6.72E+55&	1,261,274&	8,649&	81.22&	17.51\%&	UB&	10&	26.88\\
100&	5&	9.12E+63&	921,302&	1&	36.78&	24.70\%&	UB&	8&	32.90\\
100&	10&	6.50E+85&	276,122&	2,179&	21.66&	37.18\%&	UB&	10&	27.33\\
100&	20&	1.00E+98&	87,510&	1,219&	17.81&	157.94\%&	UB&	4&	31.20\\
100&	25&	1.88E+98&	61,962&	293&	17.55&	141.75\%&	UB&	1&	33.68\\
100&	50&	2.73E+78&	14,114&	938&	11.77&	0&	OPT&	0&	38.01\\
			\hline
			\multicolumn{10}{l}{$^*$ Number of feasible solutions.}\\
			\multicolumn{10}{l}{$\dagger$ Times best solution found}\\
			\multicolumn{10}{l}{$\ddagger$ Run time for all 10,000 runs for Heuristic (seconds) and one run of QAP (minutes)}\\
			\multicolumn{10}{l}{\% Percent above best found solution}\\
			\multicolumn{10}{l}{+ OPT: Optimal solution found within 5 hrs; UB: stopped with UB after 5 hrs}\\
			\multicolumn{10}{l}{$^*$ This upper bound was confirmed as optimal by CPLEX after 1,000 hours run time}\\
			
		\end{tabular}
	\end{center}
\end{table}

\begin{table}[ht!]
	\begin{center}
		\caption{\label{res} Results for variable group sizes}
		\medskip
		\begin{tabular}{|r|r|c|r|r|r||c|c|}
			\hline
&&Combin-&\multicolumn{3}{|c||}{Heuristic Runs}&\multicolumn{2}{|c|}{CPLEX}\\
\cline{4-8}			
			$n$&$k$&ations$^*$&Best&Times$\dagger$&Time$\ddagger$&\%&+\\
			\hline
12&	2&	2,047&	39,587.8&	10,000&	0.02&	0&	OPT\\
12&	3&	86,526&	19,985.8&	10,000&	0.02&	0&	OPT\\
12&	4&	611,501&	10,641.0&	6,625&	0.02&	0&	OPT\\
12&	5&	1,379,400&	7,052.4&	2,111&	0.02&	0&	OPT\\
12&	6&	1,323,652&	3,946.0&	10,000&	0.03&	0&	OPT\\
18&	2&	131,071&	108,347.8&	4,900&	0.03&	0&	OPT\\
18&	3&	64,439,010&	52,989.4&	7,908&	0.05&	0&	OPT\\
18&	4&	2.80E+09&	30,272.6&	1,620&	0.06&	0&	OPT\\
18&	5&	2.90E+10&	19,369.2&	9,596&	0.06&	0&	OPT\\
18&	6&	1.11E+11&	12,643.4&	5,301&	0.06&	0&	OPT\\
18&	7&	1.97E+11&	9,577.4&	46&	0.08&	0&	OPT\\
18&	8&	1.89E+11&	6,777.8&	581&	0.06&	0&	OPT\\
18&	9&	1.06E+11&	5,101.2&	4,976&	0.06&	0&	OPT\\
30&	2&	5.37E+08&	316,187.2&	10,000&	0.11&	0&	OPT\\
30&	3&	3.43E+13&	150,078.4&	10,000&	0.12&	0&	OPT\\
30&	5&	7.71E+18&	66,352.2&	3,436&	0.20&	0&	UB\\
30&	6&	2.99E+20&	43,550.0&	2,514&	0.20&	0&	UB\\
30&	10&	1.73E+23&	16,670.0&	1,420&	0.20&	0.67\%&	UB\\
\hline							
40&	2&	5.50E+11&	499,930&	10,000&	0.23&	0&	OPT\\
40&	4&	5.04E+22&	143,408&	9,998&	0.39&	0&	OPT\\
40&	5&	7.57E+25&	89,530&	8,166&	0.38&	0&	UB\\
40&	8&	3.17E+31&	38,576&	9,442&	0.33&	4.84\%&	UB\\
40&	10&	2.36E+33&	25,042&	467&	0.34&	5.31\%&	UB\\
50&	2&	5.63E+14&	797,668&	10,000&	0.41&	0&	OPT\\
50&	5&	7.40E+32&	165,234&	4,773&	0.64&	2.47\%&	UB\\
50&	10&	2.62E+43&	44,602&	127&	0.62&	34.76\%&	UB\\
100&	2&	6.34E+29&	3,645,284&	7,766&	2.36&	0.46\%&	UB\\
100&	4&	6.70E+58&	1,244,694&	775&	4.69&	58.01\%&	UB\\
100&	5&	6.57E+67&	850,330&	202&	5.09&	50.45\%&	UB\\
100&	10&	2.75E+93&	233,958&	1,130&	3.81&	185.89\%&	UB\\
			\hline
			\multicolumn{8}{l}{$^*$ Number of feasible solutions.}\\
			\multicolumn{8}{l}{$\dagger$ Times out of 10,000 best found}\\
			\multicolumn{8}{l}{$\ddagger$ Time in seconds for all 10,000 runs}\\
			\multicolumn{8}{l}{$^*$ Optimality verified by CPLEX}\\
			\multicolumn{8}{l}{\% Percent above best solution}\\
			\multicolumn{8}{l}{+ Optimal solution found within 5 hrs; stopped with UB after 5 hrs}\\
		\end{tabular}
	\end{center}
\end{table}

\section{\label{sec6}Computational Experiments}

We first tested the solution approaches on the \citet{WA01} instances ($n\le 30$) using Frontline System's Premium Solver in Excel. We do not report the specific results because run time cannot be measured. The problem size for the BLP formulation is limited to 18. The evolutionary algorithm was tested on   $n=12, 18, 30$ instances. The results of the evolutionary algorithm are significantly worse than the best known solution for most instances. However, if the evolutionary algorithm is repeated many times, results close to the best known or optimal solution may be obtained. As the computational results obtained by the heuristic algorithm proposed in this paper are so good, there is no reason to apply any version of Excel's solver. As is recognized by the research community, the Solver in Excel may be a good approach for instructional purposes but it is not recommended for professional use.

For testing the heuristic procedure, computer programs were coded in Fortran using double precision arithmetic and were compiled by an Intel 11.1 Fortran Compiler  using one thread with no parallel processing. They were run on a desktop with the Intel i7-6700 3.4GHz CPU processor and 16GB RAM.
The quadratic assignment problem was also solved by a Fortran program compiled by the same compiler and run on the same computer.

CPLEX was run on a virtualized Windows environment with 16 vCPUs and 128GB of vRAM.
The physical server used was a 2 CPU (8 cores each) PowerEdge R720 Intel E5-2650 CPUs with 128 GB RAM using shared storage on MD3620i via 10GB interfaces.

We experimented with the case study taken from \citet{WA01} of up to 30 cities (distances rounded to one digit after the decimal point), and for the largest U.S. cities depicted in Table \ref{cities} for $n\ge 40$ with distances calculated by (\ref{dist}) or (\ref{dist1}) rounded to the nearest integer.

Each instance was solved for both Model A and Model B by the heuristic algorithm 10,000 times in a multi-start approach incorporating GRASP. CPLEX was run only once because it does not have a random component and it finds the guaranteed optimal solution if it is not stopped prematurely. The CPLEX was stopped after 5 hours if the optimal solution was not found and the upper bound (best feasible solution found) is reported.
 
\subsection{Testing Solution Approaches for Model A}

For Model A the Quadratic Assignment Problem (QAP) was solved by the alpha-male genetic algorithm \citep{DDT18} based on the effective genetic algorithm proposed by \citet{DM13} in addition to the heuristic algorithm and CPLEX. We performed a few experiments to determine the parameters to be used for the alpha-male genetic algorithm. For the definition of parameters see  \citet{DDT18}. We found that the parameters similar to the ones used for solving the BL instances \citep{dC06} are the most effective. The population has 100 members out of which 25 are alpha males; the number of tabu search iterations is 32$n$; the number of generations is $g=5\times\max\{1000,20n\}$, and the number of population members selected for differential improvement \citep{DM13} is $R=15$. Each instance was solved 10 times by the genetic algorithm.

The results are reported in Table \ref{reseq}.
21 out of the 31 instances were solved to optimality by CPLEX in less than 5 hours of computer time. The heuristic algorithm found the best known solution in all 31 instances. The best result was obtained by the heuristic algorithm once for the $n=100, p=5$ instance and at least 293 times out of 10,000  for all other instances. In 12 instances it found the optimal solution in all 10,000 runs. Run times are very short. The largest instance was solved 10,000 times in total time of 81 seconds, i.e., 0.0081 seconds per run. The quadratic assignment solutions were quite good except for large values of $p$ when group sizes are 4 or less. It missed the best known solution in two instances. Note that CPLEX is the only approach that can guarantee that the solution found is optimal. 

\subsection{Testing Solution Approaches for Model B}

The results for Model B are depicted in Table \ref{res}.
18 out of the 30 instances were solved to optimality by CPLEX in less than 5 hours of computer time. The heuristic algorithm found the best known solution in all 30 instances. The best result was obtained by the heuristic approach at least 46 times out of 10,000 for all instances. In seven instances the optimal solution was found in all 10,000 starts. Run times by the heuristic approach are very short. The largest instance was solved 10,000 times in total time of 5 seconds, i.e., 0.0005 seconds per run. Since Model B cannot be formulated as a QAP case, the alpha-male genetic algorithm cannot be used to solve it.

\subsection{Recommendations Based on the Experiments}

The multi-start heuristic approach was extremely fast and is thus recommended as the solution method of choice unless the problem is small and a guaranteed optimal solution is sought. If reasonably long computer time is available CPLEX is recommended, and if it fails to find the optimal solution in a reasonable time (we tested some instances for two or more days and did not get optimal results), the fast heuristic can be used. The gap, i.e., the difference between the upper and lower bounds, can give an indication whether to continue to run the CPLEX. The  alpha male genetic algorithm performed well for fixed group sizes when the group sizes are not small.

\section{\label{sec7}Conclusions}

We investigated the partition of a set of items into groups when the number of groups and the distances between items are well defined. The objective function, to be minimized, is the total of the individual sums of the distances between all members of the same group. We find the partition of the set of items that minimizes the objective function. Two models are formulated and solved. In the first model the number of items in each group is given. For example, all groups must have the same number of items. In the second model there is no restriction on the number of items in each group.

We propose an optimal algorithm for each of the two problems as well as a heuristic algorithm. Problems with up to 100 items and 50 groups are tested. In the majority of instances the optimal solution was found using CPLEX. The heuristic approach, which is very fast, found the optimal solution for these cases, and it found equal or better solutions than those found by CPLEX when CPLEX was stopped after five hours. The first problem when the sizes of the groups are given, can also be formulated and solved as a quadratic assignment problem.

The heuristic algorithm performed very well. Since it is vary fast, we can also apply Tabu search \citep{GL97,G77} or simulated annealing \citep{KGV83} or other meta heuristic algorithms especially for large problems.

The heuristic algorithm can be modified to create starting solutions for many multi-facility location problems in which the demand points that are served by facilities form clusters. Seven examples for such problems are listed in the introduction.

\renewcommand{\baselinestretch}{1}
\renewcommand{\arraystretch}{1}
\large
\normalsize

\bibliographystyle{apalike}

\begin{thebibliography}{}
	
	\bibitem[Abramowitz and Stegun, 1972]{AS72}
	Abramowitz, M. and Stegun, I. (1972).
	\newblock {\em Handbook of Mathematical Functions}.
	\newblock Dover Publications Inc., New York, NY.
	
	\bibitem[Aloise, 2009]{A09}
	Aloise, D. (2009).
	\newblock {\em Exact algorithms for minimum sum-of-squares clustering}.
	\newblock PhD thesis, Ecole Polytechnique, Montreal, Canada.
	\newblock {ISBN}:978-0-494-53792-3.
	
	\bibitem[Bagirov et~al., 2015]{BODX15}
	Bagirov, A.~M., Ordin, B., Ozturk, G., and Xavier, A.~E. (2015).
	\newblock An incremental clustering algorithm based on hyperbolic smoothing.
	\newblock {\em Computational Optimization and Applications}, 61:219--241.
	
	\bibitem[Baron et~al., 2007]{BBKW07}
	Baron, O., Berman, O., Krass, D., and Wang, Q. (2007).
	\newblock The equitable location problem on the plane.
	\newblock {\em European Journal of Operational Research}, 183:578--590.
	
	\bibitem[Berman et~al., 2009]{BDTW09}
	Berman, O., Drezner, Z., Tamir, A., and Wesolowsky, G.~O. (2009).
	\newblock Optimal location with equitable loads.
	\newblock {\em Annals of Operations Research}, 167:307--325.
	
	\bibitem[Brimberg et~al., 2000]{BHMT00}
	Brimberg, J., Hansen, P., Mladenovi{\'c}, N., and Taillard, E. (2000).
	\newblock Improvements and comparison of heuristics for solving the
	uncapacitated multisource {W}eber problem.
	\newblock {\em Operations Research}, 48:444--460.
	
	\bibitem[Calik et~al., 2015]{CLY15}
	Calik, H., Labb{\'e}, M., and Yaman, H. (2015).
	\newblock p-center problems.
	\newblock In {\em Location Science}, pages 79--92. Springer.
	
	\bibitem[Church, 2019]{C19}
	Church, R.~L. (2019).
	\newblock Understanding the {W}eber location paradigm.
	\newblock In Eiselt, H.~A. and Marianov, V., editors, {\em Contributions to
		Location Analysis - In Honor of Zvi Drezner's 75th Birthday}, pages 69--88.
	Springer.
	
	\bibitem[Church and Garfinkel, 1978]{CG78}
	Church, R.~L. and Garfinkel, R.~S. (1978).
	\newblock Locating an obnoxious facility on a network.
	\newblock {\em Transportation Science}, 12:107--118.
	
	\bibitem[Daskin, 1995]{Das95}
	Daskin, M.~S. (1995).
	\newblock {\em Network and Discrete Location: Models, Algorithms, and
		Applications}.
	\newblock John Wiley \& Sons, New York.
	
	\bibitem[Daskin and Maass, 2015]{DM15}
	Daskin, M.~S. and Maass, K.~L. (2015).
	\newblock The p-median problem.
	\newblock In Laporte, G., Nickel, S., and da~Gama, F.~S., editors, {\em
		Location science}, pages 21--45. Springer.
	
	\bibitem[de~Carvalho~Jr. and Rahmann, 2006]{dC06}
	de~Carvalho~Jr., S.~A. and Rahmann, S. (2006).
	\newblock {Microarray layout as a quadratic assignment problem}.
	\newblock In Huson, D., Kohlbacher, O., Lupas, A., Nieselt, K., and Zell, A.,
	editors, {\em Proceedings of the German Conference on Bioinformatics},
	volume~83, pages 11--20, Bonn. Gesellschaft f\"ur Informatik.
	
	\bibitem[Drezner et~al., 2009]{DDG09a}
	Drezner, T., Drezner, Z., and Guyse, J. (2009).
	\newblock Equitable service by a facility: Minimizing the {G}ini coefficient.
	\newblock {\em Computers \& Operations Research}, 36:3240--3246.
	
	\bibitem[Drezner et~al., 2014]{Pap257}
	Drezner, T., Drezner, Z., and Hulliger, B. (2014).
	\newblock The quintile share ratio in location analysis.
	\newblock {\em European Journal of Operational Research}, 236:166--174.
	
	\bibitem[Drezner, 1982]{Dr82}
	Drezner, Z. (1982).
	\newblock Competitive location strategies for two facilities.
	\newblock {\em Regional Science and Urban Economics}, 12:485--493.
	
	\bibitem[Drezner, 1984]{Pap25}
	Drezner, Z. (1984).
	\newblock The p-center problem - heuristic and optimal algorithms.
	\newblock {\em Journal of the Operational Research Society}, 35:741--748.
	
	\bibitem[Drezner, 2006]{Dr06}
	Drezner, Z. (2006).
	\newblock Finding a cluster of points and the grey pattern quadratic assignment
	problem.
	\newblock {\em {OR} Spectrum}, 28:417--436.
	
	\bibitem[Drezner, 2011]{Dr11}
	Drezner, Z. (2011).
	\newblock Continuous center problems.
	\newblock In Eiselt, H.~A. and Marianov, V., editors, {\em Foundations of
		Location Analysis}, pages 63--78. Springer.
	
	\bibitem[Drezner, 2015]{Dr15}
	Drezner, Z. (2015).
	\newblock The quadratic assignment problem.
	\newblock In Laporte, G., Nickel, S., and da~Gama, F.~S., editors, {\em
		Location Science}, pages 345--363. Springer, Chum, Heidelberg.
	
	\bibitem[Drezner and Drezner, 2019]{DDT18}
	Drezner, Z. and Drezner, T.~D. (2019).
	\newblock The alpha male genetic algorithm.
	\newblock {\em {IMA} Journal of Management Mathmatics}, 30:37--50.
	
	\bibitem[Drezner and Kalczynski, 2017]{DK17}
	Drezner, Z. and Kalczynski, P. (2017).
	\newblock The continuous grey pattern problem.
	\newblock {\em Journal of the Operational Research Society}, 68:469--483.
	
	\bibitem[Drezner et~al., 2002]{DKSW02}
	Drezner, Z., Klamroth, K., Sch{\"o}bel, A., and Wesolowsky, G.~O. (2002).
	\newblock The {W}eber problem.
	\newblock In Drezner, Z. and Hamacher, H.~W., editors, {\em Facility Location:
		Applications and Theory}, pages 1--36. Springer, Berlin.
	
	\bibitem[Drezner and Misevi\v{c}ius, 2013]{DM13}
	Drezner, Z. and Misevi\v{c}ius, A. (2013).
	\newblock Enhancing the performance of hybrid genetic algorithms by
	differential improvement.
	\newblock {\em Computers \& Operations Research}, 40:1038--1046.
	
	\bibitem[Drezner et~al., 2015]{DMP14}
	Drezner, Z., Misevi\v{c}ius, A., and Palubeckis, G. (2015).
	\newblock Exact algorithms for the solution of the grey pattern quadratic
	assignment problem.
	\newblock {\em Mathematical Methods of Operations Research}, 82:85--105.
	
	\bibitem[Drezner et~al., 1986]{DrThErk}
	Drezner, Z., Thisse, J.-F., and Wesolowsky, G.~O. (1986).
	\newblock The minmax-min location problem.
	\newblock {\em Journal of Regional Science}, 26:87--101.
	
	\bibitem[Drezner and Wesolowsky, 1978]{Pap3}
	Drezner, Z. and Wesolowsky, G.~O. (1978).
	\newblock Facility location on a sphere.
	\newblock {\em Journal of the Operational Research Society}, 29:997--1004.
	
	\bibitem[Drezner and Wesolowsky, 1996]{DW96}
	Drezner, Z. and Wesolowsky, G.~O. (1996).
	\newblock Obnoxious facility location in the interior of a planar network.
	\newblock {\em Journal of Regional Science}, 35:675--688.
	
	\bibitem[Eiselt and Laporte, 1995]{EL95}
	Eiselt, H.~A. and Laporte, G. (1995).
	\newblock Objectives in location problems.
	\newblock In Drezner, Z., editor, {\em Facility Location: A Survey of
		Applications and Methods}, pages 151--180. Springer, New York, NY.
	
	\bibitem[Elzinga and Hearn, 1972]{EH72}
	Elzinga, J. and Hearn, D. (1972).
	\newblock Geometrical solutions for some minimax location problems.
	\newblock {\em Transportation Science}, 6:379--394.
	
	\bibitem[Erkut, 1990]{ER90}
	Erkut, E. (1990).
	\newblock The discrete $p$-dispersion problem.
	\newblock {\em European Journal of Operational Research}, 46:48--60.
	
	\bibitem[Erkut and Neuman, 1989]{Erk:89}
	Erkut, E. and Neuman, S. (1989).
	\newblock Analytical models for locating undesirable facilities.
	\newblock {\em European Journal of Operational Research}, 40:275--291.
	
	\bibitem[Feo and Resende, 1995]{FR95}
	Feo, T.~A. and Resende, M.~G. (1995).
	\newblock Greedy randomized adaptive search procedures.
	\newblock {\em Journal of global optimization}, 6:109--133.
	
	\bibitem[Gilmore, 1962]{Gil62}
	Gilmore, P. (1962).
	\newblock Optimal and suboptimal algorithms for the quadratic assignment
	problem.
	\newblock {\em Journal of {SIAM}}, 10:305--313.
	
	\bibitem[Glover, 1977]{G77}
	Glover, F. (1977).
	\newblock Heuristics for integer programming using surrogate constraints.
	\newblock {\em Decision Sciences}, 8:156--166.
	
	\bibitem[Glover and Laguna, 1997]{GL97}
	Glover, F. and Laguna, M. (1997).
	\newblock {\em Tabu Search}.
	\newblock Kluwer Academic Publishers, Boston.
	
	\bibitem[Hakimi, 1983]{Hak83}
	Hakimi, S.~L. (1983).
	\newblock On locating new facilities in a competitive environment.
	\newblock {\em European Journal of Operational Research}, 12:29--35.
	
	\bibitem[Hartigan and Wong, 1979]{HW79}
	Hartigan, J. and Wong, M. (1979).
	\newblock Algorithm {AS} 136: A k-means clustering algorithm.
	\newblock {\em Journal of the Royal Statistical Society. Series C (Applied
		Statistics)}, 28:100--108.
	
	\bibitem[Hotelling, 1929]{H29}
	Hotelling, H. (1929).
	\newblock Stability in competition.
	\newblock {\em Economic Journal}, 39:41--57.
	
	\bibitem[Kariv and Hakimi, 1979a]{KH79cen}
	Kariv, O. and Hakimi, S.~L. (1979a).
	\newblock An algorithmic approach to network location problems. {I}: The
	$p$-centers.
	\newblock {\em {SIAM} Journal on Applied Mathematics}, 37:513--538.
	
	\bibitem[Kariv and Hakimi, 1979b]{KH79med}
	Kariv, O. and Hakimi, S.~L. (1979b).
	\newblock An algorithmic approach to network location problems. {II}: The
	$p$-medians.
	\newblock {\em {SIAM} Journal on Applied Mathematics}, 37:539--560.
	
	\bibitem[Kirkpatrick et~al., 1983]{KGV83}
	Kirkpatrick, S., Gelat, C.~D., and Vecchi, M.~P. (1983).
	\newblock Optimization by simulated annealing.
	\newblock {\em Science}, 220:671--680.
	
	\bibitem[Koopmans and Beckmann, 1957]{KB57}
	Koopmans, T.~C. and Beckmann, M.~J. (1957).
	\newblock Assignment problems and the location of economic activities.
	\newblock {\em Econometrica}, 25:53--76.
	
	\bibitem[Kuenne and Soland, 1972]{KS72}
	Kuenne, R.~E. and Soland, R.~M. (1972).
	\newblock Exact and approximate solutions to the multisource {W}eber problem.
	\newblock {\em Mathematical Programming}, 3:193--209.
	
	\bibitem[Lawler, 1963]{Law}
	Lawler, E. (1963).
	\newblock The quadratic assignment problem.
	\newblock {\em Management Science}, 9:586--599.
	
	\bibitem[Lloyd, 1982]{L57}
	Lloyd, S. (1982).
	\newblock Least squares quantization in {PCM}.
	\newblock {\em IEEE transactions on information theory}, 28:129--137.
	
	\bibitem[MacQueen, 1967]{M67a}
	MacQueen, J. (1967).
	\newblock Some methods for classification and analysis of multivariate
	observations.
	\newblock In {\em Proceedings of the fifth Berkeley symposium on mathematical
		statistics and probability}, volume~1, pages 281--297. Oakland, CA, USA.
	
	\bibitem[Maimon, 1986]{Ma86}
	Maimon, O. (1986).
	\newblock The variance equity measure in locational decision theory.
	\newblock {\em Annals of Operations Research}, 6:147--160.
	
	\bibitem[Okabe et~al., 2000]{OKSC00}
	Okabe, A., Boots, B., Sugihara, K., and Chiu, S.~N. (2000).
	\newblock {\em Spatial Tessellations: Concepts and Applications of {V}oronoi
		Diagrams}.
	\newblock Wiley Series in Probability and Statistics. John Wiley.
	
	\bibitem[Pereira et~al., 2018]{PABM19}
	Pereira, T., Aloise, D., Brimberg, J., and Mladenovi{\'c}, N. (2018).
	\newblock Review of basic local searches for solving the minimum sum-of-squares
	clustering problem.
	\newblock In {\em Open Problems in Optimization and Data Analysis}, pages
	249--270. Springer.
	
	\bibitem[Rendl, 2002]{Rendl}
	Rendl, F. (2002).
	\newblock The quadratic assignment problem.
	\newblock In Drezner, Z. and Hamacher, H., editors, {\em Facility Location:
		Applications and Theory}. Springer, Berlin.
	
	\bibitem[Stirling, 1764]{S1764}
	Stirling, J. (1764).
	\newblock {\em Methodus differentialis, sive Tractatus de summatione et
		interpolatione serierum infinitarum. Auctore Jacobo Stirling, RSS}.
	\newblock prostat apud J. Whiston \& B. White, in Fleet-street.
	
	\bibitem[Sugihara and Iri, 1992]{SI92}
	Sugihara, K. and Iri, M. (1992).
	\newblock Construction of the voronoi diagram for ``one million" generators in
	single-precision arithmetic.
	\newblock {\em Proceedings of the IEEE}, 80:1471--1484.
	
	\bibitem[Suzuki and Okabe, 1995]{SOK95}
	Suzuki, A. and Okabe, A. (1995).
	\newblock Using {V}oronoi diagrams.
	\newblock In Drezner, Z., editor, {\em Facility Location: A Survey of
		Applications and Methods}, pages 103--118. Springer, New York.
	
	\bibitem[Sylvester, 1857]{Sylv}
	Sylvester, J. (1857).
	\newblock A question in the geometry of situation.
	\newblock {\em Quarterly Journal of Mathematics}, 1:79.
	
	\bibitem[Sylvester, 1860]{Sylv2}
	Sylvester, J. (1860).
	\newblock On {P}oncelet's approximate linear valuation of {S}urd forms.
	\newblock {\em Philosophical Magazine}, 20 (Fourth series):203--222.
	
	\bibitem[Taillard, 1995]{Ta95}
	Taillard, {\'E}.~D. (1995).
	\newblock Comparison of iterative searches for the quadratic assignment
	problem.
	\newblock {\em Location Science}, 3:87--105.
	
	\bibitem[Vorono{\"\i}, 1908]{V08}
	Vorono{\"\i}, G. (1908).
	\newblock Nouvelles applications des param{\`e}tres continus {\`a} la
	th{\'e}orie des formes quadratiques. deuxi{\`e}me m{\'e}moire. recherches sur
	les parall{\'e}llo{\`e}dres primitifs.
	\newblock {\em Journal f{\"u}r die reine und angewandte Mathematik},
	134:198--287.
	
	\bibitem[{W}eber, 1909]{W09}
	{W}eber, A. (1909).
	\newblock {\em {\"U}{b}er {d}en Standort {d}er Industrien, 1. Teil: Reine
		Theorie {d}es Standortes. English Translation: on the Location of
		Industries}.
	\newblock University of {C}hicago {P}ress, Chicago, IL.
	\newblock Translation published in 1929.
	
	\bibitem[Wesolowsky, 1993]{W93}
	Wesolowsky, G.~O. (1993).
	\newblock The {W}eber problem: History and perspectives.
	\newblock {\em Location Science}, 1:5--23.
	
	\bibitem[Winston and Albright, 2016]{WA01}
	Winston, W.~L. and Albright, S.~C. (2016).
	\newblock {\em Practical Management Science}.
	\newblock Nelson Education.
	\newblock $6^{th}$ Edition.
	
\end{thebibliography}

\end{document}